\documentclass{article}

\usepackage{amsmath}
\usepackage{graphicx}
\usepackage{hyperref}
\usepackage{cite}
\usepackage{caption}
\usepackage{subcaption}
\usepackage{authblk}

\title{Fitting an Elephant with Four non-Zero Parameters}
\author[1]{Dian Jin\thanks{Corresponding author} }
\author[2]{Junze Yuan}

\affil[1]{\textit{Research Institute for Intelligent Wearable Systems, The Hong Kong Polytechinic University}}
\affil[2]{\textit{Faculty of Engineering and Natural Sciences, ATME, Tampere University}}

\date{}

\begin{document}
	
% Title Page
\maketitle

% Abstract
\begin{abstract}
	In 1953, Enrico Fermi criticized Dyson’s model by quoting Johnny von Neumann: ``With four parameters I can fit an elephant, and with five I can make him wiggle his trunk." So far, there have been several attempts to fit an elephant using four parameters, but as the problem has not been well-defined, the current methods do not completely satisfy the requirements. This paper defines the problem and presents an attempt.
\end{abstract}

% Introduction
\section{Introduction}
In 1953, Enrico Fermi criticized Dyson’s model by quoting Johnny von Neumann: ``With four parameters I can fit an elephant, and with five I can make him wiggle his trunk."\cite{dyson2004meeting}.
This quote is intended to tell Dyson that while his model may appear complex and precise, merely increasing the number of parameters to fit the data does not necessarily imply that the model has real physical significance. 

Focus on Fermi's quotation, how can we fit an elephant with four parameters? Currently, there are some studies attempting to propose solutions.

Wei \cite{wei1975least} was the first to attempt fitting an elephant using Fourier series, but it was only with more than 20 parameters that the fit became reasonably good, but still ``may not satisfy the third-grade art teacher."

\begin{figure}[h]
	\centering
	\includegraphics[width=0.6\linewidth]{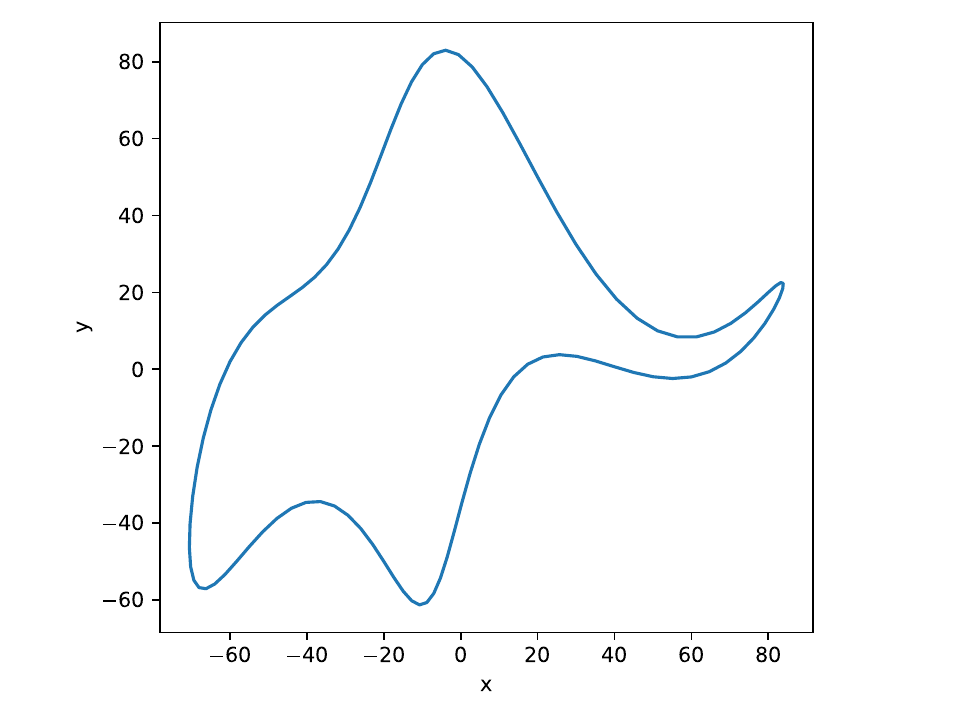}
	\caption{Elephant curve by Mayer \textit{et al.}}
	\label{fig:elephant2010}
\end{figure}

Mayer \textit{et al.} \cite{mayerDrawingElephantFour2010} used four complex numbers to fit an elephant. However, since complex numbers have both real and imaginary parts, they actually used eight parameters. Additionally, they did not count the zero parameters (if counted, it would be 20). Nevertheless, the resulting image is beautiful, with a Picasso-like style (Fig. \ref{fig:elephant2010}).

Paintadosi \cite{piantadosiOneParameterAlways2018} argues that one parameter is always enough. He constructed a function that, through a single parameter, can depict any shape. However, in essence, this work is a form of encoding, mapping the shape into a real number with precision extending to hundreds or even thousands of decimal places. For our problem, this is meaningless, although the paper's theme is that ``parameter counting" fails as a measure of model complexity.

Thus, in this paper, we define the problem and attempt to fit an elephant using four parameters under the given conditions. The code for this paper is available \footnote{https://github.com/CLaSLoVe/von-Neumann-elephant}.

% Mathematical Model
\section{Methods}

\subsection{Problem Definition}

As this problem is not strictly mathematical, some degree of ambiguity is inevitable. This definition aims to clarify the conditions under which von Neumann's elephant problem should be addressed.

Von Neumann's Elephant Problem:

Given the following requirements, fit an elephant:

\begin{enumerate}
	\item The basis function used should be widely recognized, such as the Fourier series.
	\item Ideally, the function approximation should be parameterized by exactly four parameters \( \theta = (\theta_1, \theta_2, \theta_3, \theta_4) \). Alternatively, it can be parameterized such that exactly four of these parameters \( \theta_i \) (\( i = 1, 2, 3, 4 \)) are non-zero.
	\item The resulting curve \( f(x,y; \theta) \) should approximate the shape of an elephant, exhibiting key morphological characteristics of an elephant.
\end{enumerate}

% Parameter Determination
\subsection{Model fitting}

Just as in previous research, we still utilize Fourier series, but in polar coordinates:

$$
r(\theta) = \sum_{k=1}^{\infty} \left( a_k \cos(k\theta) + b_k \sin(k\theta) \right) + c
$$
where $c$ is a constant.

Since von Neumann did not specify whether the elephant he referred to was 2D or 3D, frontal or lateral, previous research has considered the lateral view of a 2D elephant. We consider the frontal view of a 2D elephant. This approach has several advantages: First, biological organisms typically exhibit a certain degree of symmetry, which allows for the incorporation of prior knowledge during fitting, thereby reducing the number of parameters required for fitting. Second, the most memorable features of an elephant, such as tusks, trunk, and large ears, are primarily located on the face, making the frontal view a good choice.

\begin{figure}[h]
	\centering
	\includegraphics[width=0.6\linewidth]{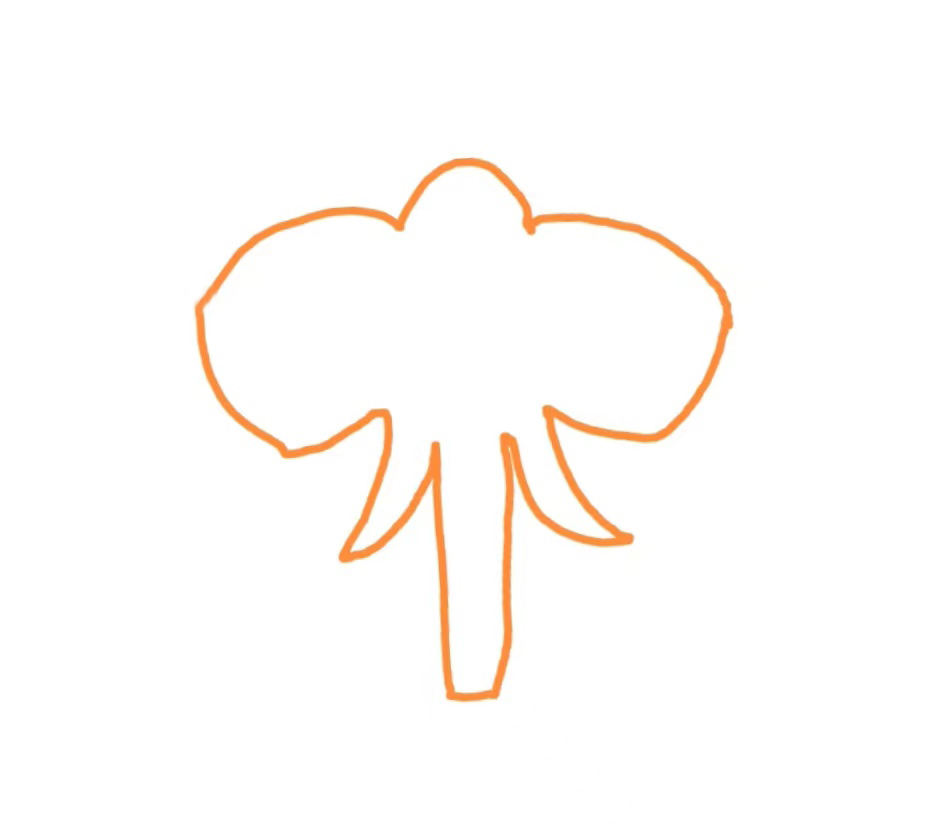}
	\caption{Abstract picture of an elephant}
	\label{fig:elephant}
\end{figure}

We drew a front-view picture of an elephant using a single-stroke drawing method, as shown in Fig. \ref{fig:elephant}.

Due to the presence of symmetry, we consider \( r(\theta) = r(-\theta) \):

\begin{align*}
	r(\theta) &= \sum_{k=1}^{\infty} \left( a_k \cos(k\theta) + b_k \sin(k\theta) \right) + c \\
	r(-\theta) &= \sum_{k=1}^{\infty} \left( a_k \cos(-k\theta) + b_k \sin(-k\theta) \right) + c
\end{align*}

So,
$$
\sum_{k=1}^{\infty} \left( a_k \cos(k\theta) + b_k \sin(k\theta) \right) = \sum_{k=1}^{\infty} \left( a_k \cos(k\theta) - b_k \sin(k\theta) \right)
$$

Thus we can obtain that $b_k = 0$, and the series simplifies to:

$$
r(\theta) = \sum_{k=1}^{\infty} a_k \cos(k\theta) + c
$$

During the fitting process, we uniformly sampled $ n $ points ($ n = 100 $) from Fig. \ref{fig:elephant}, with each point represented in polar coordinates $(r_0, \theta_0)$.

As for $ c $, we estimated it using the mean value $\bar{r_0}$ of all the sampled points.

We used the least squares method as the loss function:

$$
L(a_k) = \sum_{i=1}^{n} \left( r_i - f(\theta_i; a_k) \right)^2
$$

\begin{figure}[h]
	\centering
	\begin{subfigure}[b]{0.49\linewidth}
		\centering
		\includegraphics[width=\linewidth]{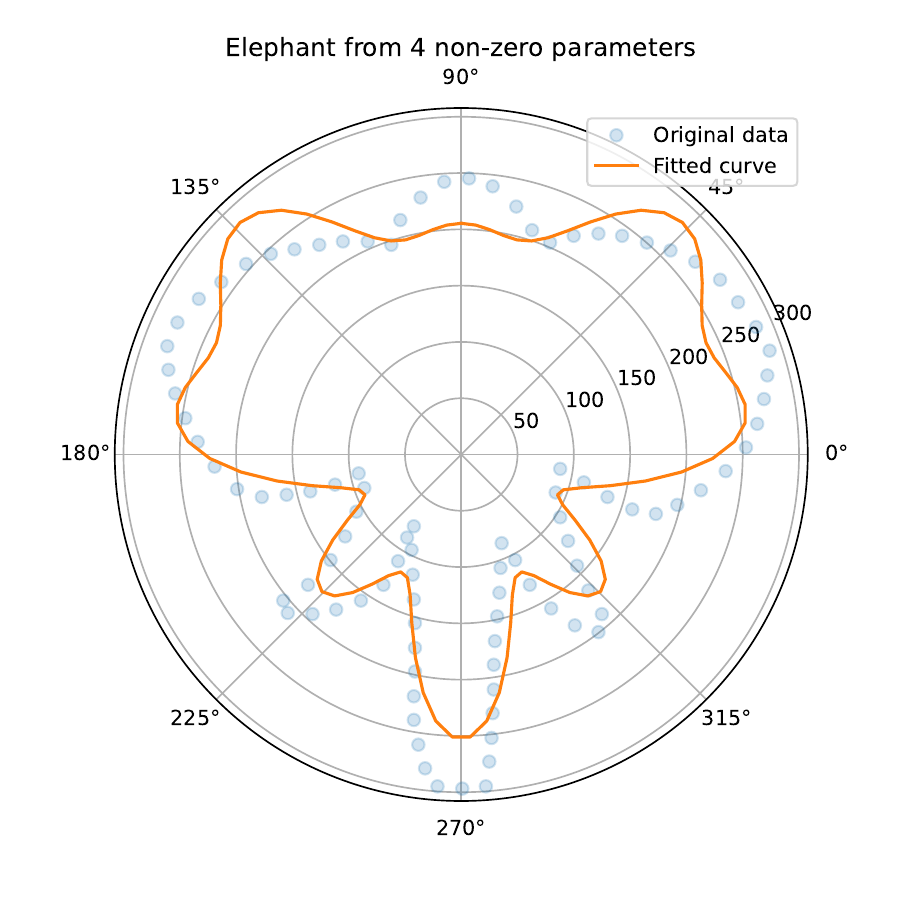}
		\caption{Result of fitting elephant curve}
	\end{subfigure}
	\hfill
	\begin{subfigure}[b]{0.49\linewidth}
		\centering
		\includegraphics[width=\linewidth]{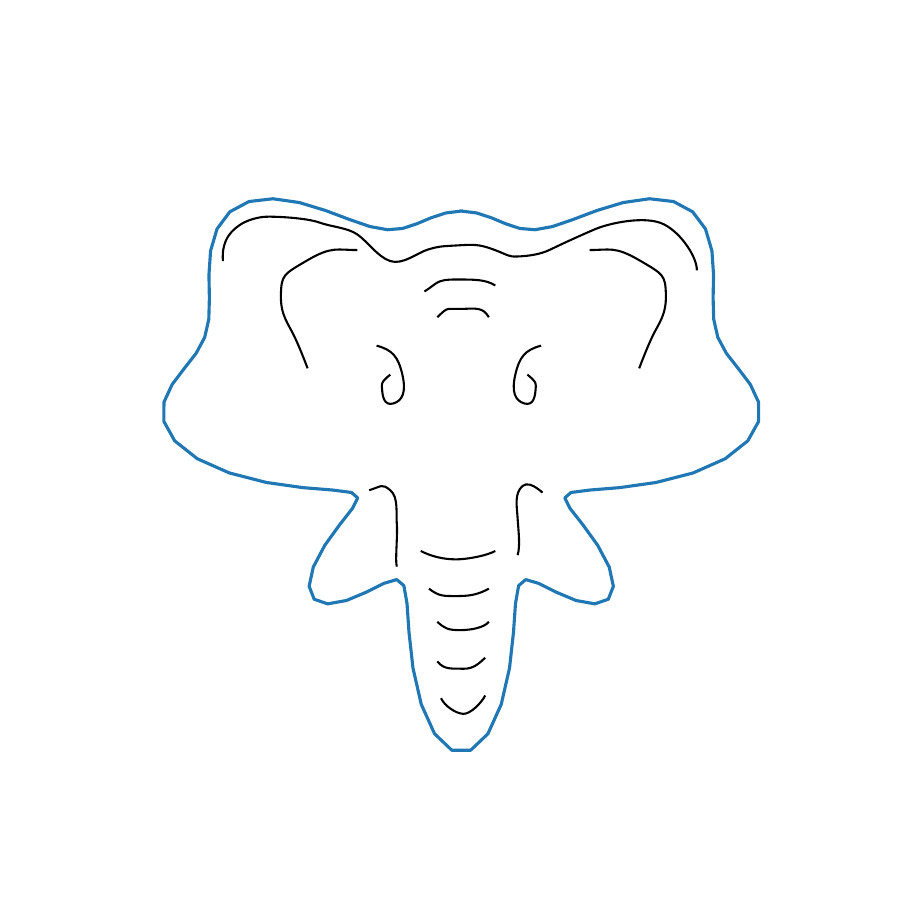}
		\caption{The front view of an imagined elephant}
	\end{subfigure}
	\caption{Fitting result}
	\label{fig:result}
\end{figure}

% Results
\section{Results}

The result of our fitted curve is shown in Eq. \ref{eq:result}. The fitted curve is shown in Fig. \ref{fig:result}(a). We used four non-zero parameters (or eight parameters) to fit the front view of an elephant. The front view of an imagined elephant is shown  Fig. \ref{fig:result}(b).

\begin{equation}
	\begin{aligned}
		r = & 47.84 \cdot \cos(1 \cdot \theta) 
		-51.12 \cdot \cos(3 \cdot \theta) \\
		-&20.43 \cdot \cos(7 \cdot \theta)
		+31.58 \cdot \cos(8 \cdot \theta) \\
		+& \bar{r_0}
	\end{aligned}
	\label{eq:result}
\end{equation}

% Discussion
\section{Discussion}
In this paper, we define the problem of von Neumann's elephant and fit the front view of an elephant using four non-zero parameters. Without the defined constraints, this problem would become trivial and meaningless; for example, you could construct an encoder to encode any number into a pattern, or you could create an ``elephant'' coordinate system or define an ``elephant'' curve using zero parameters. 

There are several limitations of this paper:

\begin{enumerate}
	\item It only satisfies a weaker condition, i.e., using four non-zero parameters instead of four parameters.
	\item It only shows the front view of the elephant, and the presented content is not as comprehensive as in previous studies.
\end{enumerate}

\section*{Acknowledgement}

We would like to thank Jianyong Jin for the illustration of the elephant. We also thank ``Tai Le Mao Alice" for introducing this topic to us.

% References
\bibliographystyle{IEEEtran}
\bibliography{lib}

% Generated by IEEEtran.bst, version: 1.14 (2015/08/26)
\begin{thebibliography}{1}
\providecommand{\url}[1]{#1}
\csname url@samestyle\endcsname
\providecommand{\newblock}{\relax}
\providecommand{\bibinfo}[2]{#2}
\providecommand{\BIBentrySTDinterwordspacing}{\spaceskip=0pt\relax}
\providecommand{\BIBentryALTinterwordstretchfactor}{4}
\providecommand{\BIBentryALTinterwordspacing}{\spaceskip=\fontdimen2\font plus
\BIBentryALTinterwordstretchfactor\fontdimen3\font minus
  \fontdimen4\font\relax}
\providecommand{\BIBforeignlanguage}[2]{{%
\expandafter\ifx\csname l@#1\endcsname\relax
\typeout{** WARNING: IEEEtran.bst: No hyphenation pattern has been}%
\typeout{** loaded for the language `#1'. Using the pattern for}%
\typeout{** the default language instead.}%
\else
\language=\csname l@#1\endcsname
\fi
#2}}
\providecommand{\BIBdecl}{\relax}
\BIBdecl

\bibitem{dyson2004meeting}
F.~Dyson \emph{et~al.}, ``A meeting with enrico fermi,'' \emph{Nature}, vol.
  427, no. 6972, pp. 297--297, 2004.

\bibitem{wei1975least}
J.~Wei, ``Least square fitting of an elephant,'' \emph{Chemtech}, vol.~5,
  no.~2, pp. 128--129, 1975.

\bibitem{mayerDrawingElephantFour2010}
J.~Mayer, K.~Khairy, and J.~Howard, ``Drawing an elephant with four complex
  parameters,'' \emph{American Journal of Physics}, vol.~78, no.~6, pp.
  648--649, Jun. 2010.

\bibitem{piantadosiOneParameterAlways2018}
S.~T. Piantadosi, ``One parameter is always enough,'' \emph{AIP Advances},
  vol.~8, no.~9, p. 095118, Sep. 2018.

\end{thebibliography}
	
\end{document}